\documentclass{amsart}

\usepackage{graphicx,amssymb,mathrsfs,amsmath,color,fancyhdr}
\usepackage[all]{xy}
\newtheorem{theorem}{Theorem}[section]

\theoremstyle{definition}

\theoremstyle{remark}

\numberwithin{equation}{section}

\begin{document}

\title[The elementary symmetric functions of a reciprocal polynomial sequence]
{The elementary symmetric functions of a reciprocal polynomial sequence}

\begin{abstract}
Erd\"{o}s and Niven proved in 1946 that for any positive integers $m$ and $d$,
there are at most finitely many integers $n$ for which at least one of the elementary
symmetric functions of $1/m, 1/(m+d), ..., 1/(m+(n-1)d)$ are integers. Recently, Wang
and Hong refined this result by showing that if $n\geq 4$, then none of the elementary
symmetric functions of $1/m, 1/(m+d), ..., 1/(m+(n-1)d)$ is an integer for any positive
integers $m$ and $d$. Let $f$ be a polynomial of degree at least $2$ and of
nonnegative integer coefficients. In this paper, we show that none of the elementary
symmetric functions of $1/f(1), 1/f(2), ..., 1/f(n)$ is an integer except for
$f(x)=x^{m}$ with $m\geq2$ being an integer and $n=1$.
\end{abstract}
\author{Yuanyuan Luo}
\address{Mathematical College, Sichuan University, Chengdu 610064, P.R. China}
\email{yuanyuanluoluo@163.com}
\author{Shaofang Hong$^*$}
%    Address of record for the research reported here
\address{Mathematical College, Sichuan University, Chengdu 610064, P.R. China}
%\curraddr{}
\email{sfhong@scu.edu.cn; s-f.hong@tom.com; hongsf02@yahoo.com}
\author{Guoyou Qian}
%    Address of record for the research reported here
\address{Center for Combinatorics, Nankai University, Tianjin 300071, P.R. China}
\email{qiangy1230@163.com; qiangy1230@gmail.com}
\author{Chunlin Wang}
\address{Mathematical College, Sichuan University, Chengdu 610064, P.R. China}
\email{wdychl@126.com}
\thanks{$^*$S. Hong is the corresponding author and was supported
partially by National Science Foundation of
China Grant \#11371260. G. Qian was supported partially by Postdoctoral
Science Foundation of China Grant \#2013M530109}
\maketitle

\section{Introduction}
Let $n$ be a positive integer and $f(x)$ be a polynomial of integer coefficients
such that $f(m)\ne 0$ for any integer $m\ge 1$. For any integer $k$ with $1\le k\le n$,
we denote by $\sigma_{k, f}(n)$ the $k$-th elementary symmetric functions of
$1/f(1), 1/f(2), ..., 1/f(n)$. That is,
$$\sigma_{k, f}(n):=\sum\limits_{1\leq i_{1}< i_{2}<\cdots<i_{k}\leq n}
\prod\limits_{j=1}^{k}\frac{1}{f(i_{j})}.$$
A well-known result says that if $n\geq 2$ and $f(x)=x$, then the harmonic sum
$\sigma_{1, f}(n)$ cannot be an integer. More generally, if $n\geq 2$ and $f(x)=ax+b$
with $a$ and $b$ being positive integers, then the sum $\sigma_{1, f}(n)$ is not an integer.
In 1946, Erd\"{o}s and Niven \cite{[EN]} extended this result by showing that if $f(x)=ax+b$ with
$a$ and $b$ being positive integers, then there are at most finitely many integers
$n$ for which at least one element in the set $S(f, n):=\{\sigma_{1, f}(n), \sigma_{2, f}(n), ...,
\sigma_{n, f}(n)\}$ is an integer. In 2012, Chen and Tang \cite{[CT]} proved that each element
of $S(f, n)$ is not an integer if $f(x)=x$ and $n\geq 4$. Wang and Hong \cite{[WH1]}
showed that none of the elements in $S(f, n)$ is an integer if $f(x)=2x-1$ and $n\geq 2$.
Recently, Wang and Hong \cite{[WH2]} refined the theorem of Erd\"{o}s and Niven \cite{[EN]}
by showing that if $f(x)=ax+b$ with $a$ and $b$ being positive integers and $n\geq 4$,
then all the elements in $S(f, n)$ are not integers. An interesting problem
naturally arises: Does the similar result hold when $f(x)$ is a polynomial of
nonnegative integer coefficients and of degree at least two?

In this paper, our main goal is to answer the above problem. In fact, we determine
all the finite progressions $\{f(i)\}_{i=1}^{n}$ with $f(x)$ being of nonnegative
coefficients such that one or more elements
in $S(f, n)$ are integers. In other words, we have the following result.
\begin{theorem}\label{Thm}
Let $f$ be a polynomial of nonnegative integer coefficients and of degree at least two.
Let $n$ and $k$ be integers such that $1\le k\le n$. Then $\sigma_{k, f}(n)$ is not an
integer except for the case $f(x)=x^{m}$ with $m\geq2$ being an integer and $k=n=1$,
in which case, $\sigma_{k, f}(n)$ is an integer.
\end{theorem}

Evidently, Theorem {\ref{Thm}} answers completely the above problem.
In the next section, we will give the proof of Theorem {\ref{Thm}}.
A conjecture is proposed in the last section.

\section{Proof of Theorem \ref{Thm}}
This section is devoted to the proof of Theorem {\ref{Thm}}. To do so,
we first list two known identities about the values of Riemann zeta
function at 2 and 4 (see, for example, \cite{[K]}):
$$\zeta(2)=\sum\limits_{j=1}^{\infty}\frac{1}{j^{2}}=\frac{\pi^{2}}{6} \ {\rm and} \
\zeta(4)=\sum\limits_{j=1}^{\infty}\frac{1}{j^{4}}=\frac{\pi^{4}}{90}.$$
Then we can easily see that $1<\zeta(2)<2$. Notice that $\sigma_{k, f}(n)>0$
for any integer $n\geq 1$.

We can now give the proof of Theorem \ref{Thm}.\\
\\
{\it Proof of Theorem \ref{Thm}.}
Let $f(x)=a_m x^m +a_{m-1}x^{m-1}+...+a_0$ with $a_m \geq 1$ and $m\geq 2$ being integers.
First we let $k\geq 2$. It follows from the hypotheses $a_m \geq 1$ and $m\geq 2$ that
$f(r)\geq r^2$ for any positive integer $r$. Since $\zeta(2)<2$, we deduce that
\begin{align}\label{1}
\sigma_{k+1, f}(n)&=\sum\limits_{1\leq i_{1}<\cdots<i_{k+1}\leq n}\prod
\limits_{j=1}^{k+1}\frac{1}{f(i_{j})}\\
&=\sum\limits_{1\leq i_{1}<\cdots<i_{k}\leq n-1}\Big(\prod\limits_{j=1}^{k}
\frac{1}{f(i_{j})}\Big)\Big(\sum\limits_{i_{k+1}=i_{k}+1}^{n}\frac{1}{f(i_{k+1})}\Big)\nonumber\\
&\leq\sum\limits_{1\leq i_{1}<\cdots<i_{k}\leq n-1}\Big(\prod\limits_{j=1}^{k}
\frac{1}{f(i_{j})}\Big)\Big(\sum\limits_{i_{k+1}=2}^{\infty}\frac{1}{i_{k+1}^{2}}\Big)\nonumber\\
&=(\zeta(2)-1)\sum\limits_{1\leq i_{1}<\cdots<i_{k}\leq n-1}\prod\limits_{j=1}^{k}
\frac{1}{f(i_{j})}\nonumber\\
&=(\zeta(2)-1)\sigma_{k, f}(n-1)\nonumber\\
&<\sigma_{k, f}(n-1)<\sigma_{k, f}(n)\nonumber.
\end{align}
So for any given integer $n$, $\sigma_{k, f}(n)$ is decreasing as $k$ increases.
On the other hand, we have
\begin{align}\label{2}
\sigma_{2, f}(n)&=\sum\limits_{1\leq i_{1}<i_{2}\leq n}\frac{1}{f(i_{1})f(i_{2})}\\
&\leq\sum\limits_{1\leq i_{1}<i_{2}\leq n}\frac{1}{i_{1}^2i_{2}^{2}}\nonumber
\end{align}
\begin{align}
&<\sum\limits_{i_{2}>i_{1}\ge 1}\frac{1}{i_{1}^2i_{2}^{2}}\nonumber\\
&=\frac{1}{2}\Big(\Big(\sum\limits_{j=1}^{\infty}\frac{1}{j^{2}}\Big)^{2}-
\sum\limits_{j=1}^{\infty}\frac{1}{j^{4}}\Big)\nonumber\\
&=\frac{1}{2}(\zeta(2)^2 -\zeta(4))=\frac{\pi^{4}}{120}<1\nonumber.
\end{align}
Thus, by (\ref{1}) and (\ref{2}), we obtain that $0<\sigma_{k, f}(n)<1$ if $2\leq k\leq n$.
This concludes that $\sigma_{k, f}(n)$ is not an integer if $k\geq 2$. So Theorem {\ref{Thm}}
is true for the case that $k\geq 2$.

In what follows we let $k=1$. First we assume that $f$ contains only one term, namely
$f(x)=ax^{m}$, where $m\geq2$ and $a\ge 1$. Clearly, if $a\ge2$, then
$$0<\sigma_{1, f}(n)\le\frac{1}{2}\sum_{j=1}^{n}\frac{1}{j^2}
<\frac{1}{2}\sum_{j=1}^{\infty }\frac{1}{j^2}=\frac{\pi^2}{12}<1.$$
If $a=1$, then $f(x)=x^{m}$. It follows that $\sigma_{1, f}(1)=1$ and
$$1<\sigma_{1, f}(n)\leq\sum\limits_{j=1}^{\infty}\frac{1}{j^{2}}=\frac{\pi^{2}}{6}<2$$
for any integer $n\ge 2$. Hence for any $n\ge 2$, $\sigma_{1, f}(n)$ is not
an integer if $f(x)=a_mx^{m}$ with $m\geq2$ and $a_m\ge 1$.

Now we suppose that $f$ contains at least two terms. Then one may let
$f(x)=a_m x^m +a_{m-1}x^{m-1}+...+a_0,$ where $m\ge 2, a_m\ge 1$ and
$\max(a_0,..., a_{m-1})\ge 1$. We divide the proof into the following three cases:

{\bf Case 1.} $m=2, a_1=0, a_0=a_2=1$. Then $f(x)=x^2 +1$. By a simple
calculation we see that $\sigma_{1, f}(12)<1$, $\sigma_{1, f}(13)>1$.
So we can conclude that $0<\sigma_{1, f}(n)\le \sigma_{1, f}(12)<1$ if $n\le 12$, and
$$1<\sigma_{1, f}(13)\le \sigma_{1, f}(n)<\sum\limits_{j=1}^{\infty}
\frac{1}{j^{2}+1}<\sum\limits_{j=1}^{\infty}\frac{1}{j^{2}}=\zeta(2)<2$$
if $n\ge 13.$ Thus $\sigma_{1, f}(n)$ is not an integer in this case.

{\bf Case 2.} $m=2$, $a_1=0$ and $\max (a_0, a_2)\ge 2$.
Then for any positive integer $j$, one can deduce that $f(j)=a_2j^2+a_0\geq j^2+2$. It then follows that
$$0<\sigma_{1, f}(1)\leq \frac{1}{3}<1, 0<\sigma_{1, f}(2)\leq \frac{1}{3}+\frac{1}{6}<1$$
and
$$0<\sigma_{1, f}(n)\leq \sum\limits_{j=1}^{n}\frac{1}{j^{2}+2}<\frac{1}{3}+\frac{1}{6}+
\sum\limits_{j=3}^{n}\frac{1}{(j-1)j}=\frac{1}{3}+\frac{1}{6}+\frac{1}{2}-\frac{1}{n}<1$$
if $n\geq 3.$ Namely, $\sigma_{1, f}(n)$ is not an integer in this case.

{\bf Case 3.} Either $m=2$ and $a_1\ge 1$, or $m\ge 3$. If $m\ge 3$, since $f(x)$
contains at least two terms, it follows that there is an integer $l$
with $0\le l<m$ such that $a_l\ge 1$. Hence for any positive integer $j$,
we derive that
$$f(j)\ge a_mj^m+a_lj^l\ge j^3+1\ge j^2+j$$
if $m\ge 3$. If $m=2$ and $a_1\ge 1$, then for any positive integer $j$, we have
$f(j)=a_2 j^2+a_1j+a_0\geq j^2+j$. Based on the above discussions, we can deduce that
$$0<\sigma_{1, f}(n)=\sum\limits_{j=1}^{n}\frac{1}{f(j)}\leq
\sum\limits_{j=1}^{n}\frac{1}{j^2+j}=1-\frac{1}{n+1}<1.$$
So $\sigma_{1, f}(n)$ is not an integer in this case.

This completes the proof of Theorem {\ref{Thm}} for the case that $k=1$.
So Theorem {\ref{Thm}} is proved. \hfill$\Box$

\section{Remarks}
In this section, we raise the following conjecture as the conclusion of this paper.\\
\\
{\bf Conjecture 3.1.}
{\it Let $f(x)$ be a polynomial of integer coefficients such that $f(m)\ne 0$
for any positive integer $m$. Then there is a positive
integer $N$ such that for any integer $n\ge N$ and for all integers $k$
with $1\le k\le n$, $\sigma_{k, f}(n)$ is not an integer.}\\
\\
Clearly, by \cite{[EN]} (or \cite{[WH2]}) and Theorem 1.1 we know that
Conjecture 3.1 is true if $f(x)$ is of nonnegative integer coefficients.
Further, by \cite{[WH2]} one can derive that Conjecture 3.1 holds if $f(x)=ax-b$,
where $a$ and $b$ are integers such that $a>b>0$. But it is kept open
for the case that either $f(x)=ax-b$ with $a$ and $b$ being integers such that $0<a<b$,
or $f(x)$ is of degree greater than 2 and contains negative coefficients
but its leading coefficient is positive.
\bibliographystyle{amsplain}

\end{document}